\documentclass{emspublic}

\usepackage{epsfig}
\usepackage{amsmath,amssymb,amsfonts,amsthm,amscd}

\DeclareSymbolFont{EulerScript}{U}{eus}{m}{n}
\DeclareSymbolFontAlphabet\mathscr{EulerScript}

\newcommand{\N}{\mathscr N}
\newcommand{\Z}{\mathbb Z}
\newcommand{\Q}{\mathbb Q}
\newcommand{\R}{\mathbb R}
\newcommand{\inte}{\hbox{int}\,}
\newcommand{\lk}{\ell k}
\newcommand{\im}{\hbox{Im}\,}
\renewcommand{\ker}{\hbox{Ker}\,}

\begin{document}
	
\title{The Conway function of a splice}
\author{David Cimasoni\thanks{Supported by the Swiss National Science Foundation.}}
\shortauthor{David Cimasoni}
\affiliation{Section de Math\'ematiques, Universit\'e de Gen\`eve, 2--4 rue du Li\`evre, 1211 Gen\`eve 24, Switzerland}
\received{}

\maketitle

\begin{abstract}
We give a closed formula for the Conway function of a splice in terms of the Conway function of its splice components.
As corollaries, we refine and generalize results of Seifert, Torres, and Sumners-Woods.
\end{abstract}

\amsprimary{57M25}

\section{Introduction}

The connected sum, the disjoint sum, and cabling are well-known operations on links. As pointed out by Eisenbud and Neumann
\cite{E-N}, these are special cases of an operation which they call {\em splicing\/}. Informally, the splice of two links
$L'$ and $L''$ along components $K'\subset L'$ and $K''\subset L''$ is the link
$(L'\setminus K')\cup(L''\setminus K'')$ obtained by pasting the exterior of $K'$
and the exterior of $K''$ along their boundary torus (see Section \ref{section:prel} for a precise definition).
But splicing is not only a natural generalization of classical operations. Indeed, Eisenbud and Neumann gave the following
reinterpretation of the Jaco-Shalen and Johannson splitting theorem: Any irreducible link in an integral homology sphere can be
expressed as the result of splicing together a collection of Seifert links and hyperbolic links, and the minimal way of
doing this is unique (see \cite[Theorem 2.2]{E-N}).

Given such a natural operation, it is legitimate to ask how invariants of links behave under splicing.
Eisenbud and Neumann gave the answer for several invariants, including the multivariable Alexander polynomial
(see \cite[Theorem 5.3]{E-N}). For an oriented ordered link $L$ with $n$ components in an integral homology sphere,
this invariant is an element
of the ring $\Z[t^{\pm 1}_1,\dots,t^{\pm 1}_n]$, well-defined up to multiplication by
$\pm t^{\nu_1}_1\cdots t^{\nu_n}_n$ for integers $\nu_1,\dots,\nu_n$. Now, there exists a refinement of the Alexander
polynomial called the {\em Conway function\/}, which is a well-defined rational function $\nabla_L\in\Z(t_1,\dots,t_n)$.
This invariant was first introduced by Conway in \cite{Con}
and formally defined by Hartley \cite{Har} for links in $S^3$. The extension to links in any integral homology sphere
is due to Turaev \cite{Tu1}.

In this paper, we give a closed formula for the Conway function of a splice $L$ in terms of the Conway function of its
splice components $L'$ and $L''$. This result can be considered as a refinement of \cite[Theorem 5.3]{E-N}.
As applications, we refine well-known formulae of Seifert, Torres, and Sumners-Woods.

The paper is organized as follows. In Section \ref{section:prel}, we define the splicing and the Conway function 
as a refined torsion. Section \ref{section:res} contains the statement of the main result (Theorem \ref{thm})
and a discussion of several of its consequences
(Corollaries \ref{cor:1} to \ref{cor:4}). Finally, Section \ref{section:pr} deals with the proof of Theorem \ref{thm}.

\section{Preliminaries}
\label{section:prel}

In this section, we begin by recalling the definition of the splicing operation as introduced in \cite{E-N}.
Then, we define the torsion of a chain complex and the sign-determined torsion of a homologically oriented 
$CW$-complex following \cite{Tu3}. We finally recall Turaev's definition of the Conway function, refering to \cite{Tu1}
for further details.

\bigskip

\noindent{\bf Splice.} Let $K$ be a knot in a $\Z$-homology sphere $\Sigma$ and let $\N(K)$ be
a closed tubular neighborhood of $K$ in $\Sigma$. A pair $\mu,\lambda$
of oriented simple closed curves in $\partial\N(K)$ is said to be a {\em standard meridian and longitude\/} for $K$ if $\mu\sim 0$,
$\lambda\sim K$ in $H_1(\N(K))$ and $\lk_\Sigma(\mu,K)=1$, $\lk_\Sigma(\lambda,K)=0$,
where $\lk_\Sigma(-,-)$ is the linking number in $\Sigma$. Note that this pair is unique up to isotopy.

Consider two oriented links
$L'$ and $L''$ in $\Z$-homology spheres $\Sigma'$ and $\Sigma''$, and choose components $K'$ of $L'$ and $K''$ of $L''$. Let
$\mu',\lambda'\subset\partial\N(K')$ and $\mu'',\lambda''\subset\partial\N(K'')$ be standard meridians and longitudes. Set
$$
\Sigma=(\Sigma'\setminus\inte\N(K'))\cup(\Sigma''\setminus\inte\N(K'')),
$$
where the pasting homeomorphism maps $\mu'$ onto $\lambda''$ and $\lambda'$ onto $\mu''$.
The link $(L'\setminus K')\cup(L''\setminus K'')$ in $\Sigma$ is called the {\em splice of $L'$ and $L''$ along $K'$ and 
$K''$\/}. The manifold $\Sigma$ is easily seen to be a $\Z$-homology sphere. However, even if $\Sigma'=\Sigma''=S^3$, $\Sigma$
might not be the standard sphere $S^3$. This is the reason for considering links in $\Z$-homology spheres from the start.

Let us mention the following easy fact (see \cite[Proposition 1.2]{E-N} for the proof).
\begin{lemma}\label{lemma:lk}
Given any component $K_i$ of $L'\setminus K'$ and $K_j$ of $L''\setminus K''$,
$$
\lk_{\Sigma}(K_i,K_j)=\lk_{\Sigma'}(K',K_i)\,\lk_{\Sigma''}(K'',K_j).
$$
\end{lemma}

\bigskip

\noindent{\bf Torsion of chain complexes.}
Given two bases $c,c'$ of a finite-dimensional vector space on a field $F$, let $[c/c']\in F^*$ be the determinant
of the matrix expressing the vectors of the basis $c$ as linear combination of vectors in $c'$.

Let $C=(C_m\to C_{m-1}\to\dots\to C_0)$ be a finite-dimensional chain complex over a field $F$, such that for
$i=0,\dots,m$, both $C_i$ and $H_i(C)$ have a distinguished basis. Set
$$
\beta_i(C)=\sum_{r\le i}\dim H_r(C),\quad\gamma_i(C)=\sum_{r\le i}\dim C_r.
$$
Let $c_i$ be the given basis of $C_i$ and $h_i$ a
sequence of vectors in $\ker(\partial_{i-1}\colon C_i\to C_{i-1})$ whose projections in $H_i(C)$ form the given basis
of $H_i(C)$. Let $b_i$ be a sequence of vectors in $C_i$ such that $\partial_{i-1}(b_i)$ forms a basis of
$\im(\partial_{i-1})$. Clearly, the sequence $\partial_i(b_{i+1})h_ib_i$ is a basis of $C_i$.
The {\em torsion of the chain complex $C$} is defined as
$$
\tau(C)=(-1)^{\vert C\vert}\prod_{i=0}^m[\partial_i(b_{i+1})h_ib_i/c_i]^{(-1)^{i+1}}\in F^*,
$$
where $\vert C\vert=\sum_{i=0}^m\beta_i(C)\gamma_i(C)$. It turns out that $\tau(C)$ depends on
the choice of bases of $C_i$, $H_i(C)$, but does not depend on the choice of $h_i$, $b_i$.

We shall need the following lemma, which follows easily from \cite[Lemma 3.4.2]{Tu1} and \cite[Remark 1.4.1]{Tu2}.

\begin{lemma}\label{mult}
Let $0\to C'\to C\to C''\to 0$ be an exact sequence of finite-dimensional chain complexes of length $m$ over $F$. Assume
that the vector spaces $C_i',C_i,C''_i$ and $H_i(C'),H_i(C),H_i(C'')$ have distinguished bases. Then,
$$
\tau(C)=(-1)^{\mu+\nu}\,\tau(C')\,\tau(C'')\,\tau(\mathscr{H})\prod_{i=0}^m[c_i'c_i''/c_i]^{(-1)^{i+1}},
$$
where $\mathscr{H}$ is the based acyclic chain complex
$$
\mathscr{H}=(H_m(C')\to H_m(C)\to H_m(C'')\to\dots\to H_0(C')\to H_0(C)\to H_0(C''))
$$
and
\begin{eqnarray*}
\nu&\!\!=\!\!&\sum_{i=0}^m\gamma_i(C'')\gamma_{i-1}(C'),\\
\mu&\!\!=\!\!&\sum_{i=0}^m\left(\beta_i(C)+1\right)\left(\beta_i(C')+\beta_i(C'')\right)+\beta_{i-1}(C')\beta_i(C'').
\end{eqnarray*}
\end{lemma}

\bigskip

\noindent{\bf Sign-determined torsions of $\mathbf{CW}$-complexes.}
Consider now a finite $CW$-complex $X$ and a ring homomorphism $\varphi\colon\Z[H]\to F$, where $H=H_1(X;\Z)$. Assume that
$X$ is {\em homologically oriented\/}, that is, is endowed with a preferred orientation $\omega$ of the real vector space
$H_*(X;\R)=\bigoplus_{i\ge 0}H_i(X;\R)$. To this triple $(X,\varphi,\omega)$, we associate a {\em sign-determined
torsion\/} $\tau^{\varphi}(X,\omega)\in F/\varphi(H)$ as follows. Consider the maximal abelian covering $\widehat X\to X$ and endow
$\widehat X$ with the induced $CW$-structure. Clearly, the cellular chain complex $C(\widehat X)$ is a complex of
$\Z[H]$-modules and $\Z[H]$-linear homomorphisms. Viewing $F$ as a $\Z[H]$-module via the homomorphism $\varphi$, one has
the chain complex over $F$
$$
C^\varphi(X)=F\otimes_{\Z[H]}C(\widehat X).
$$
If this complex is not acyclic, set $\tau^{\varphi}(X,\omega)=0$. Assume $C^\varphi(X)$ is acyclic. Choose a family
$\hat e$ of cells of $\widehat X$ such that over each cell of $X$ lies exactly one cell of $\hat e$. Orient and order these
cells in an arbitrary way. This yields a basis of $C(\widehat X)$ over $\Z[H]$, and thus a basis of $C^\varphi(X)$
over $F$ and a torsion $\tau(C^\varphi(X))\in F^*$. Moreover, the orientation and the order of the cells of
$\hat e$ induce an orientation and an order for the cells of $X$, and thus a basis for the cellular chain complex
$C(X;\R)$. Choose a basis $h_i$ of $H_i(X;\R)$ such that the basis $h_0 h_1\dots h_{\dim X}$ of $H_*(X;\R)$ is
positively oriented with respect to $\omega$. Consider the torsion $\tau(C(X;\R))\in\R^*$ of the
resulting based chain complex with based homology. Denote by $\tau_0$ its sign and set
$$
\tau^{\varphi}(X,\omega)=\tau_0\,\tau(C^\varphi(X))\in F^*.
$$
It turns out that $\tau^{\varphi}(X,\omega)$ only depends on $(X,\varphi,\omega)$ and $\hat e$. Furthermore, its
class in $F/\varphi(H)$ does not depend on $\hat e$. This class is the sign-determined torsion of $X$.

\bigskip

\noindent{\bf The Conway function.}
Let $L=K_1\cup\dots\cup K_n$ be an oriented link in an oriented integral homology sphere $\Sigma$. Let $\N_i$ be a closed
tubular neighborhood of $K_i$ for $i=1,\dots,n$, and let $X$ be a cellular structure on
$\Sigma\setminus\sqcup_{i=1}^n\inte\N_i$.
Recall that $H=H_1(X;\Z)$ is a free abelian group on $n$ generators $t_1,\dots,t_n$
represented by the meridians of $K_1,\dots,K_n$. Let $F=Q(H)$ be the field of fractions of 
the ring $\Z[H]$, and let $\varphi\colon\Z[H]\hookrightarrow Q(H)$ be the standard
inclusion. Finally, let $\omega$ be the homology orientation of $X$ given by the basis of $H_*(X;\R)$
$$
\left([pt],t_1,\dots,t_n,[\partial\N_1],\dots,[\partial\N_{n-1}]\right),
$$
where $\partial\N_i$ is oriented as the boundary of $\N_i$. (The space $\N_i$ inherits the orientation of $\Sigma$.)
Consider the sign-determined torsion $\tau^{\varphi}(X,\omega)\in Q(H)/H$. It turns out to satisfy the equation
$$
\tau^{\varphi}(X,\omega)(t_1^{-1},\dots,t_n^{-1})=
(-1)^n\,t_1^{\nu_1}\cdots t_n^{\nu_n}\,\tau^{\varphi}(X,\omega)(t_1,\dots,t_n)
$$
for some integers $\nu_1,\dots,\nu_n$. The {\em Conway function of the link $L$\/} is the rational function
$$
\nabla_L(t_1,\dots,t_n)=-t_1^{\nu_1}\cdots t_n^{\nu_n}\,\tau^{\varphi}(X,\omega)(t_1^2,\dots,t_n^2)\in Q(H).
$$
Note that it satisfies the equation $\nabla_L(t_1^{-1},\dots,t_n^{-1})=(-1)^n\nabla_L(t_1,\dots,t_n)$.
The {\em reduced Conway function of $L$\/} is the one-variable Laurent polynomial
$$
\Omega_L(t)=(t-t^{-1})\nabla_L(t,\dots,t)\in\Z[t^{\pm 1}].
$$

We shall need one basic property of $\nabla_L$ known as the {\em Torres formula\/} (see e.g. \cite{B-L} for a proof).
\begin{lemma}\label{Torres}
Let $L=K_1\cup\dots\cup K_n$ be an oriented link in an integral homology sphere $\Sigma$.
If $L'$ is obtained from $L$ by removing the component $K_1$, then
$$
\nabla_L(1,t_2,\dots,t_n)=
(t_2^{\ell_2}\cdots t_n^{\ell_n}-t_2^{-\ell_2}\cdots t_n^{-\ell_n})
\nabla_{L'}(t_2,\dots,t_n),
$$
where $\ell_i=\lk_\Sigma(K_1,K_i)$ for $2\le i\le n$.
\end{lemma}

\section{The results}
\label{section:res}

\begin{theorem}\label{thm}
Let $L=K_1\cup\dots\cup K_n$ be the splice of $L'=K'\cup K_1\cup\dots\cup K_m$ and $L''=K''\cup K_{m+1}\cup\dots\cup K_n$
along $K'$ and $K''$, with $n>m\ge 0$. Let $\ell_i'$ and $\ell_j''$ denote the linking numbers $\lk_{\Sigma'}(K',K_i)$
and $\lk_{\Sigma''}(K'',K_j)$. Then,
$$
\nabla_L(t_1,\dots,t_n)=\nabla_{L'}(t_{m+1}^{\ell''_{m+1}}\cdots t^{\ell''_n}_n,t_1,\dots,t_m)\nabla_{L''}(t_1^{\ell'_1}\cdots t^{\ell'_m}_m,t_{m+1},\dots,t_n),
$$
unless $m=0$ and $\ell_1''=\dots=\ell_n''=0$, in which case
$$
\nabla_L(t_1,\dots,t_n)=\nabla_{L''\setminus K''}(t_1,\dots,t_n).
$$
\end{theorem}

Let us give several corollaries of this result, starting with the following {\em Seifert-Torres formula\/}
for the Conway function. (The corresponding formula for the Alexander polynomial was proved by Seifert \cite{Sei}
in the case of knots, and by Torres \cite{Tor} for links.)

\begin{corollary}\label{cor:1}
Let $K$ be a knot in a $\Z$-homology sphere, and $\N(K)$ a closed tubular neighborhood of $K$. Consider an
orientation-preserving homeomorphism $f\colon\N(K)\to S^1\times D^2$ that maps $K$ onto $S^1\times\{0\}$ and a standard
longitude onto $S^1\times\{1\}$. If $L=K_1\cup\dots\cup K_n$ is a link in the interior of $\N(K)$ with
$K_i\sim \ell_i\cdot K$ in $H_1(\N(K))$, then
$$
\nabla_L(t_1,\dots,t_n)=\Omega_K(t_1^{\ell_1}\cdots t_n^{\ell_n})\nabla_{f(L)}(t_1,\dots,t_n).
$$
\end{corollary}
\begin{proof}
The link $L$ is nothing but the splice of $K$ and $\mu\cup f(L)$ along $K$ and $\mu$, where $\mu$ denotes
a meridian of $S^1\times D^2$. If $\ell_i\neq 0$ for some $1\le i\le n$, Theorem \ref{thm} and Lemma \ref{Torres} give
\begin{eqnarray*}
\nabla_L(t_1,\dots,t_n)&=&\nabla_K(t_1^{\ell_1}\cdots t_n^{\ell_n})\nabla_{\mu\cup f(L)}(1,t_1,\dots,t_n)\\
	&=&\nabla_K(t_1^{\ell_1}\cdots t_n^{\ell_n})(t_1^{\ell_1}\cdots t_n^{\ell_n}-t_1^{-\ell_1}\cdots t_n^{-\ell_n})
	\nabla_{f(L)}(t_1,\dots,t_n)\\
	&=&\Omega_K(t_1^{\ell_1}\cdots t_n^{\ell_n})\nabla_{f(L)}(t_1,\dots,t_n).
\end{eqnarray*}
On the other hand, if $\ell_i=0$ for all $i$, then Theorem \ref{thm} implies
$$
\nabla_L(t_1,\dots,t_n)=\nabla_{f(L)}(t_1,\dots,t_n).
$$
Since $K$ is a knot, $\Omega_K(1)=1$ and the corollary is proved.
\end{proof}

\begin{corollary}\label{cor:2}
Assuming the notation of Corollary \ref{cor:1}, we have
$$
\Omega_L(t)=\Omega_K(t^\ell)\Omega_{f(L)}(t)
$$
if $L\sim \ell\cdot K$ in $H_1(\N(K))$.
\end{corollary}
\begin{proof}
Use Corollary \ref{cor:1} and the definition of $\Omega_L(t)$.
\end{proof}

Let $p,q$ be coprime integers. Recall that a {\em $(p,q)$-cable\/} of a knot $K$ is a knot on $\partial\N(K)$
homologous to $p\cdot\lambda+q\cdot\mu$, where $\mu,\lambda$ is a standard meridian and longitude for $K$.
We have the following refinement and generalization of \cite[Theorems 5.1 to 5.4]{S-W}.

\begin{corollary}\label{cor:3}
Let $L=K_1\cup\dots\cup K_n$ be an oriented link in a $\Z$-homology sphere, and let $\ell_i=\lk_\Sigma(K_i,K_n)$
for $1\le i\le n-1$.
Consider the link $L'$ obtained from $L$ by adding $d$ parallel copies of a $(p,q)$-cable of $K_n$. Then,
$$
\nabla_{L'}(t_1,\dots,t_{n+d})=\left(t_n^qT^p-t_n^{-q}T^{-p}\right)^d
\nabla_L(t_1,\dots,t_{n-1},t_n(t_{n+1}\cdots t_{n+d})^p)
$$
and
$$
\nabla_{L'\setminus K_n}(t_1,\dots,\widehat t_n,\dots,t_{n+d})=\frac{\left(T^p-T^{-p}\right)^d}{T-T^{-1}}
\nabla_L(t_1,\dots,t_{n-1},(t_{n+1}\cdots t_{n+d})^p),
$$
where $T=t_1^{\ell_1}\cdots t_{n-1}^{\ell_{n-1}}(t_{n+1}\cdots t_{n+d})^q $.
\end{corollary}
\begin{proof}
Let $L''$ be the link in $S^3$ consisting of $d$ parallel copies $K_{n+1}\cup\dots\cup K_{n+d}$ of a $(p,q)$-torus knot on
a torus $Z$, together with the oriented cores $K''$, $K'$ of the two solid tori bounded by $Z$.
Let us say that $K''$ is the core such that $\lk(K'',K_{n+i})=p$, and that $K'$ satisfies 
$\lk(K',K_{n+i}')=q$ for $1\le i\le d$. By \cite{Cim}, the Conway function of $L''$ is given by
$$
\nabla_{L''}(t'',t',t_{n+1},\dots,t_{n+d})=
\left({t''}^p{t'}^q(t_{n+1}\cdots t_{n+d})^{pq}-{t''}^{-p}{t'}^{-q}(t_{n+1}\cdots t_{n+d})^{-pq}\right)^d.
$$
The link $L'$ is the splice of $L$ and $L''$ along $K_n$ and $K''$. By Theorem \ref{thm},
$$
\nabla_{L'}(t_1,\dots,t_{n+d})=\nabla_L(t_1,\dots,t_{n-1},t_n(t_{n+1}\cdots t_{n+d})^p)
\nabla_{L''}(t_1^{\ell_1}\cdots t_n^{\ell_n},t_n,t_{n+1},\dots,t_{n+d})
$$
leading to the first result. The value of $\nabla_{L\setminus K'}$ then follows from Lemma \ref{Torres}.
\end{proof}

Note that for $d=1$ and $(p,q)=(2,1)$, the second equality of Corollary \ref{cor:3} is nothing but Turaev's
`Doubling Axiom' (see \cite[p. 154]{Tu1} and \cite[p. 105]{Tu2}).

\begin{corollary}\label{cor:4}
If $L$ is the connected sum of $L'=K'_1\cup K_2\cup\dots\cup K_m$ and $L''=K''_1\cup K_{m+1}\cup\dots\cup K_n$ along
$K'_1$ and $K''_1$, then
$$
\nabla_L(t_1,\dots,t_n)=(t_1-t_1^{-1})\nabla_{L'}(t_1,\dots,t_m)\nabla_{L''}(t_1,t_{m+1},\dots,t_n).
$$
\end{corollary}
\begin{figure}[Htb]
   \begin{center}
     \epsfig{figure=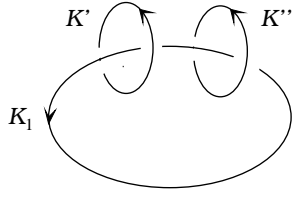,height=2.5cm}
     \caption{\footnotesize The link $\widetilde L$ in the proof of Corollary \ref{cor:4}.}
     \label{fig}
   \end{center}
\end{figure}
\begin{proof}
Consider the link $\widetilde L=K'\cup K''\cup K_1$ illustrated in Figure \ref{fig}. The link $L$ can be understood as
the splice of $\widetilde L$ and $L'$ along $K'$ and $K'_1$, itself spliced with $L''$ along $K''$ and $K''_1$. Since
$\nabla_{\widetilde L}(t',t'',t_1)=t_1-t_1^{-1}$, the result
follows easily from Theorem \ref{thm}.
\end{proof}

Finally, note that Theorem \ref{thm} can also be understood as a generalization of the Torres formula
(Lemma \ref{Torres}). Indeed, $L'$ is the splice of $L$ and the trivial knot $K$ along $K_1$ and $K$. Theorem
\ref{thm} then leads to the Torres formula. Nevertheless, it should not be considered as a corollary of our result,
since we will make use of this formula in our proof.

\section{Proof of Theorem \ref{thm}}
\label{section:pr}

The first step of the proof consists of reducing the general case to a simpler situation using Lemma \ref{Torres}.

\begin{lemma}
Assume that Theorem \ref{thm} holds when $\ell'_i\neq 0$ for some $1\le i\le m$ and $\ell''_j\neq 0$
for some $m+1\le j\le n$. Then Theorem \ref{thm} always holds.
\end{lemma}
\begin{proof}
Let us first assume that $m>0$.
Set $\widetilde L'=K_0'\cup L'$, where $K_0'$ is an oriented knot in $\Sigma'\setminus L'$ such that
$\lk_{\Sigma'}(K_0',K')=\lk_{\Sigma'}(K_0',K_1)=1$ and $\lk_{\Sigma'}(K_0',K_i)=0$ for
$i>1$. Similarly, set $\widetilde L''=K_0''\cup L''$, where $K_0''$ is an oriented knot in $\Sigma''\setminus
L''$ such that $\lk_{\Sigma''}(K_0'',K'')=\lk_{\Sigma''}(K_0'',K_{m+1})=1$ and $\lk_{\Sigma''}(K_0'',K_j)=0$ for
$j>m+1$ (recall that $n>m$). Consider the splice $\widetilde L$ of $\widetilde L'$ and $\widetilde L''$ along $K'$
and $K''$. This splice satisfies the conditions of the statement, so Theorem \ref{thm} can be applied, giving
$$
\nabla_{\widetilde L}(t_0',t_0'',t_1,\dots,t_n)=\nabla_{\widetilde L'}(t_0',t_0''T',t_1,\dots,t_m)
\nabla_{\widetilde L''}(t_0'',t_0'T'',t_{m+1},\dots,t_n),
$$
where $T'=t_{m+1}^{\ell''_{m+1}}\cdots t^{\ell''_n}_n$ and $T''=t_1^{\ell'_1}\cdots t_m^{\ell'_m}$.
Setting $t_0'=1$ and applying Lemmas \ref{Torres} and \ref{lemma:lk}, we get that
$(t_0''t_1T'-(t_0''t_1T')^{-1})\nabla_{\widetilde L\setminus K_0'}(t_0'',t_1,\dots,t_n)$ is equal to the product
$$
(t_0''t_1T'-(t_0''t_1T')^{-1})\nabla_{L'}(t_0''T',t_1,\dots,t_m)\nabla_{\widetilde L''}(t_0'',T'',t_{m+1},\dots,t_n).
$$
Since $t_0''t_1T'-(t_0''t_1T')^{-1}\neq 0$, the equation can be divided by this factor. Setting $t_0''=1$, we see that
$(t_{m+1}T''-(t_{m+1}T'')^{-1})\nabla_L(t_1,\dots,t_n)$ is equal to
$$
(t_{m+1}T''-(t_{m+1}T'')^{-1})\nabla_{L'}(T',t_1,\dots,t_m)\nabla_{L''}(T'',t_{m+1},\dots,t_n).
$$
Since $t_{m+1}T''-(t_{m+1}T'')^{-1}\neq 0$, the case $m>0$ is proved.

Assume now that $m=0$ and $\ell''_j\neq 0$ for some $1\le j\le n$.
Set $\widetilde L'=K_0'\cup K'$, where $K_0'$ is a meridian of $K'$. Consider
the splice $\widetilde L$ of $\widetilde L'$ and $L''$ along $K'$
and $K''$. Since $\widetilde L'$ is not a knot and the case $m>0$ holds, we can apply Theorem \ref{thm}. This gives
$$
\nabla_{\widetilde L}(t_0',t_1,\dots,t_n)=\nabla_{\widetilde L'}(t_0',T')\nabla_{L''}(t_0',t_1,\dots,t_n),
$$
where $T'=t_1^{\ell''_{1}}\cdots t^{\ell''_n}_n$. Setting $t_0'=1$, we get
$$
(T'-{T'}^{-1})\nabla_L(t_1,\dots,t_n)=(T'-{T'}^{-1})\nabla_{L'}(T')\nabla_{L''}(1,t_1,\dots,t_n).
$$
Since $T'-{T'}^{-1}\neq 0$, the case $m=0$ is settled if $\ell''_j\neq 0$ for some $1\le j\le n$.

Finally, assume that $m=0$ and $\ell''_1=\dots=\ell''_n=0$. 
Set $\widetilde L''=K_0''\cup L''$, where $K_0''$ is a meridian of $K''$. Since $\lk_{\Sigma''}(K_0'',K'')\neq 0$,
the theorem can be applied to the splice $\widetilde L$ of $L'$ and $\widetilde L''$ along $K'$
and $K''$:
\begin{eqnarray*}
\nabla_{\widetilde L}(t_0'',t_1,\dots,t_n)&=&\nabla_{L'}(t_0'')\nabla_{\widetilde L''}(t_0',1,t_1,\dots,t_n)\\
	&=&\nabla_{L'}(t_0'')(t_0''-t_0''^{-1})\nabla_{\widetilde L''\setminus K''}(t_0'',t_1,\dots,t_n)\\
	&=&\Omega_{L'}(t_0'')\nabla_{\widetilde L''\setminus K''}(t_0'',t_1,\dots,t_n).
\end{eqnarray*}
Setting $t_0''=1$ and using the fact that $\Omega_{L'}(1)=1$, it follows
$$
(t_1-t_1^{-1})\nabla_L(t_1,\dots,t_n)=(t_1-t_1^{-1})\nabla_{L''\setminus K''}(t_1,\dots,t_n)
$$
and the lemma is proved.
\end{proof}

So, let us assume that $\ell'_i\neq 0$ for some $1\le i\le m$ and $\ell''_j\neq 0$ for some $m+1\le j\le n$.

Let $X$ be a cellular decomposition of $\Sigma\setminus\inte\N(L)$ having $X'=\Sigma'\setminus\inte\N(L')$,
$X''=\Sigma''\setminus\inte\N(L'')$ and $T=\partial\N(K')=\partial\N(K'')$ as subcomplexes. Note that $H=H_1(X;\Z)$ is
free abelian with basis $t_1,\dots,t_n$ represented by meridians of $K_1,\dots,K_n$. Similarly, $H'=H_1(X';\Z)$ has basis
$t',t_1,\dots,t_m$, $H''=H_1(X'';\Z)$ has basis $t'',t_{m+1},\dots,t_n$ and $H_T=H_1(T;\Z)$ has basis $t',t''$.
Moreover, the inclusion homomorphism $H'\to H$ is given by $t_i\mapsto t_i$ for $1\le i\le m$ and
$t'\mapsto t_{m+1}^{\ell''_{m+1}}\cdots t^{\ell''_n}_n$. Since $\ell''_j\neq 0$ for some $m+1\le j\le n$, it is injective.
Therefore, it induces a monomomorphism $j'\colon\Z[H']\to \Z[H]$ which fits in the commutative diagram
$$
\begin{CD}
\Z[H']@>{j'}>>\Z[H]\\
@V{\varphi'}VV @V{\varphi}VV\\
Q(H')@>{i'}>>Q(H),
\end{CD}
$$
where $\varphi$ (resp. $\varphi'$) denotes the standard inclusion of $\Z[H]$ (resp. $\Z[H']$) into its field of fractions. 
Similarly, the inclusion homomorphisms $H''\to H$ and $H_T\to H$ are injective, inducing
$$
\begin{CD}
\Z[H'']@>{j''}>>\Z[H]\\
@V{\varphi''}VV @V{\varphi}VV\\
Q(H')@>{i''}>>Q(H)
\end{CD}
\qquad\hbox{and}\qquad
\begin{CD}
\Z[H_T]@>{j_T}>>\Z[H]\\
@V{\varphi_T}VV @V{\varphi}VV\\
Q(H')@>{i_T}>>Q(H).
\end{CD}
$$

Let $\N_i=\N(K_i)$, $\N'=\N(K')$ and $\N''=\N(K'')$ be closed tubular neighborhoods.
Let $\omega$, $\omega'$, $\omega''$ be the homology orientations of $X$, $X'$, $X''$ given by the basis
\begin{eqnarray*}
h_0h_1h_2&=&\left([pt],t_1,\dots,t_n,[\partial\N_1],\dots,[\partial\N_{n-1}]\right),\\
h'_0h'_1h'_2&=&\left([pt],t',t_1,\dots,t_m,[\partial\N'],[\partial\N_1],\dots,[\partial\N_{m-1}]\right),\\
h''_0h''_1h''_2&=&\left([pt],t'',t_{m+1},\dots,t_n,[\partial\N''],[\partial\N_{m+1}],\dots,[\partial\N_{n-1}]\right)
\end{eqnarray*}
of $H_*(X;\R)$, $H_*(X';\R)$, $H_*(X'';\R)$, respectively. Finally, let $\omega_T$ be the homology orientation of $T$
given by the basis $h_0^Th_1^Th_2^T=([pt],t',t'',[\partial\N'])$ of $H_*(T;\R)$.

\begin{lemma}\label{torsion}
$$
\tau^{\varphi}(X,\omega)\,i_T(\tau^{\varphi_T}(T,\omega_T))=
i'(\tau^{\varphi'}(X',\omega'))\,i''(\tau^{\varphi''}(X'',\omega''))\in Q(H)/H.
$$
\end{lemma}
\begin{proof}
Let $p\colon\widehat X\to X$ be the universal abelian covering of $X$. Endow $\widehat X$ with the induced cellular
structure. We have the exact sequence of cellular chain complexes over $\Z[H]$
$$
0\to C(p^{-1}(T))\to C(p^{-1}(X'))\oplus C(p^{-1}(X''))\to C(\widehat X)\to 0.
$$
If $\widehat X'\to X'$ is the universal abelian covering of $X'$, then $C(p^{-1}(X'))=\Z[H]\otimes_{\Z[H']}C(\widehat X')$,
where $\Z[H]$ is a $\Z[H']$-module via the homomorphism $j'$. Therefore,
\begin{eqnarray*}
Q(H)\otimes_{\Z[H]}C(p^{-1}(X'))&=&Q(H)\otimes_{\Z[H]}\left(\Z[H]\otimes_{\Z[H']}C(\widehat X')\right)\\
	&=&Q(H)\otimes_{\Z[H']}C(\widehat X')\\
	&=&C^{\varphi\circ j'}(X')=C^{i'\circ\varphi'}(X').
\end{eqnarray*}
Similarly, we have $Q(H)\otimes_{\Z[H]}C(p^{-1}(X''))=C^{i''\circ\varphi''}(X'')$ and
$Q(H)\otimes_{\Z[H]}C(p^{-1}(T))=C^{i_T\circ\varphi_T}(T)$.
This gives the exact sequence of chain complexes over $Q(H)$
$$
0\to C^{i_T\circ\varphi_T}(T)\to C^{i'\circ\varphi'}(X')\oplus C^{i''\circ\varphi''}(X'')
\to C^{\varphi}(X)\to 0.
$$

Since the inclusion homomorphism $H_T\to H$ is non-trivial, the complex $C^{i_T\circ\varphi_T}(T)$ is acyclic (see the
proof of \cite[Lemma 1.3.3]{Tu1}). By the long exact sequence associated with the sequence of complexes given above, $C^{\varphi}(X)$
is acyclic if and only if $C^{i'\circ\varphi'}(X')$ and $C^{i''\circ\varphi''}(X'')$ are acyclic. Clearly, this is
equivalent to asking that $C^{\varphi'}(X')$ and $C^{\varphi''}(X'')$ are acyclic. Therefore,
$\tau^{\varphi_T}(T,\omega_T)\neq 0$ and
$$
\tau^{\varphi}(X,\omega)=0\quad\Longleftrightarrow\quad\tau^{\varphi'}(X',\omega')=0\;\hbox{ or }\;
\tau^{\varphi''}(X',\omega'')=0.
$$
Hence, the lemma holds in this case, and it may be assumed that $C^{\varphi}(X)$, $C^{i'\circ\varphi'}(X')$ and $C^{i''\circ\varphi''}(X'')$ are acyclic.

Choose a family $\hat e$ of cells of $\widehat X$ such that over each cell of $X$ lies exactly one cell of $\hat e$.
Orient these cells in an arbitrary way, and order them by counting first the cells over $T$, then the cells
over $X'\setminus T$, and finally the cells over $X''\setminus T$.
This yields $Q(H)$-bases $\hat c$, $\hat c^T$, $\hat c'$, $\hat c''$ for
$C^{\varphi}(X)$, $C^{i_T\circ\varphi_T}(T)$, $C^{i'\circ\varphi'}(X')$, $C^{i''\circ\varphi''}(X'')$, and $\R$-bases $c$,
$c^T$, $c'$, $c''$ for $C(X;\R)$, $C(T;\R)$, $C(X';\R)$ and $C(X'';\R)$. Applying Lemma \ref{mult} to the exact
sequence of based chain complexes above, we get
$$
\tau(C^{i'\circ\varphi'}(X')\oplus C^{i''\circ\varphi''}(X''))=
(-1)^{\nu(T,X)}\,\tau(C^{i_T\circ\varphi_T}(T))\tau(C^{\varphi}(X))\,(-1)^\sigma,
$$
where $\nu(T,X)=\sum_i\gamma_i(C^{i_T\circ\varphi_T}(T))\gamma_{i-1}(C^{\varphi}(X))=
\sum_i\gamma_i(C(T))\gamma_{i-1}(C(X))$ and
$\sigma=\sum_i(\#\hat c_i'-\#\hat c_i^T)\#\hat c_i^T=\sum_i(\# c_i'-\# c_i^T)\# c_i^T$.
Using Lemma  \ref{mult} and the exact sequence
$$
0\to C^{i'\circ\varphi'}(X')\to C^{i'\circ\varphi'}(X')\oplus C^{i''\circ\varphi''}(X'')\to C^{i''\circ\varphi''}(X'')\to 0,
$$
we get
$$
\tau(C^{i'\circ\varphi'}(X')\oplus C^{i''\circ\varphi''}(X''))=
(-1)^{\nu(X',X'')}\tau(C^{i'\circ\varphi'}(X'))\tau(C^{i''\circ\varphi''}(X'')).
$$
Therefore,
$$
\tau(C^{\varphi}(X))\tau(C^{i_T\circ\varphi_T}(T))=
(-1)^{N}\tau(C^{i'\circ\varphi'}(X'))\tau(C^{i''\circ\varphi''}(X'')),
$$
where $N=\nu(T,X)+\nu(X',X'')+\sigma$. By functoriality of the torsion (see e.g. \cite[Proposition 3.6]{Tu2}), $\tau(C^{i'\circ\varphi'}(X'))=i'(\tau(C^{\varphi'}(X')))$, $\tau(C^{i''\circ\varphi'}(X''))=i''(\tau(C^{\varphi''}(X'')))$
and $\tau(C^{i_T\circ\varphi_T}(T))=i_T(\tau(C^{\varphi_T}(T)))$. Hence,
$$
\tau(C^{\varphi}(X))i_T(\tau(C^{\varphi_T}(T)))=
(-1)^{N}i'(\tau(C^{\varphi'}(X')))i''(\tau(C^{\varphi''}(X''))).\eqno{(\star)}
$$

Now, consider the exact sequences
\begin{eqnarray*}
&0\to C(T;\R)\to C(X';\R)\oplus C(X'';\R)\to C(X;\R)\to 0&\;\hbox{and}\\
&0\to C(X';\R)\to C(X';\R)\oplus C(X'';\R)\to C(X'';\R)\to 0,&
\end{eqnarray*}
and set $\beta_i(-)=\beta_i(C(-\,;\R))$. Lemma \ref{mult} gives the equations
\begin{eqnarray*}
\tau(C(X';\R)\oplus C(X'';\R))&\!\!=\!\!&(-1)^{\mu+\nu(T,X)}\tau(\mathscr{H})\tau(C(T;\R))\tau(C(X;\R))(-1)^\sigma,\\
\tau(C(X';\R)\oplus C(X'';\R))&\!\!=\!\!&(-1)^{\tilde\mu+\nu(X',X'')}\tau(C(X';\R))\tau(C(X'';\R)),
\end{eqnarray*}
where
\begin{eqnarray*}
\mu&\!\!=\!\!&\sum_i\left(\beta_i(X')+\beta_i(X'')+1\right)\left(\beta_i(T)+\beta_i(X)\right)+\beta_{i-1}(T)\beta_i(X),\\
\tilde\mu&\!\!=\!\!&\sum_i\left(\beta_i(X')+\beta_i(X'')+1\right)\left(\beta_i(X')+\beta_i(X'')\right)+
\beta_{i-1}(X')\beta_i(X''),\\
\end{eqnarray*}
and $\mathscr{H}$ is the based acyclic complex
$$
\mathscr{H}=(H_2(T;\R)\to\dots\to H_0(T;\R)\to H_0(X';\R)\oplus H_0(X'';\R)\to H_0(X;\R)).
$$
Therefore,
$$
\tau(C(X;\R))\tau(C(T;\R))=(-1)^{M}\tau(\mathscr{H})\tau(C(X';\R))\tau(C(X'';\R)),\eqno{(\star\star)}
$$
where $M=\mu+\tilde\mu+\nu(T,X)+\nu(X',X'')+\sigma$.
By equations $(\star)$ and $(\star\star)$,
$$
\tau^{\varphi}(X,\omega)\,i_T(\tau^{\varphi_T}(T,\omega_T))=(-1)^{\mu+\tilde\mu}sign(\tau(\mathscr{H}))
i'(\tau^{\varphi'}(X',\omega'))\,i''(\tau^{\varphi''}(X'',\omega''))
$$
in $Q(H)/H$, and we are left with the proof that $sign(\tau(\mathscr{H}))=(-1)^{\mu+\tilde\mu}$.
Since $\beta_i(T)+\beta_i(X)+\beta_i(X')+\beta_i(X'')$ is even for all $i$, as well as $\beta_i(T)$ and $\beta_i(X'')$
for $i\ge 2$, we have
\begin{eqnarray*}
\mu+\tilde\mu&\equiv&\textstyle{\sum_i}\,\beta_{i-1}(T)\beta_i(X)+\beta_{i-1}(X')\beta_i(X'')\pmod{2}\\
	&\equiv&\beta_0(T)\beta_1(X)+\beta_0(X')\beta_1(X'')\pmod{2}\\
	&\equiv&m+1\pmod{2}.
\end{eqnarray*}
Furthermore, the acyclic complex $\mathscr{H}$ splits into three short exact sequences
$$
0\to H_i(T;\R)\stackrel{f_i}{\to} H_i(X';\R)\oplus H_i(X'';\R)\stackrel{g_i}{\to}H_i(X;\R)\to 0,
$$
for $i=0,1,2$. Therefore, $\tau(\mathscr{H})=\prod_{i=0}^2[f_i(h^T_i)r_i(h_i)/h_i'h_i'']^{(-1)^i}$, where
$r_i$ satisfies $g_i\circ r_i=id$.  We have $f_0(h^T_0)r_0(h_0)=([pt]\oplus-[pt],[pt]\oplus 0)$ and
$h_0'h_0''=([pt]\oplus 0,0\oplus [pt])$. Hence,
$$
[f_0(h^T_0)r_0(h_0)/h_0'h_0'']=
\begin{vmatrix}
\phantom{-}1&1\cr
-1&0
\end{vmatrix}
=1.
$$
Furthermore,
$$
[f_1(h^T_1)r_1(h_1)/h_1'h_1'']=
\begin{vmatrix}
1\cr &\ell_1'&1\cr &\vdots& &\ddots\cr &\ell_m'& & &1\cr &-1\cr -\ell''_{m+1}& & & & &1\cr \vdots& & & & & &\ddots\cr
-\ell_n''& & & & & & &1
\end{vmatrix}
=(-1)^{m+1}.
$$
Finally, using the equality $[\partial N']+ [\partial N_1]+\dots+[\partial N_m]=0$ in $H_2(X';\R)$, we have
$$
[f_2(h^T_2)r_2(h_2)/h_2'h_2'']=
\begin{vmatrix}
1& & & &-1\cr &1& & &-1\cr & &\ddots& &\vdots\cr & & &1&-1\cr 1\cr & & & & &1\cr & & & & & &\ddots\cr & & & & & & &1
\end{vmatrix}
=1.
$$
So $\tau(\mathscr H)=(-1)^{m+1}$ and the lemma is proved.
\end{proof}

It is easy to show that $\tau^{\varphi_T}(T,\omega_T)=\pm 1\in Q(H_T)/H_T$ (see \cite[Lemma 1.3.3]{Tu1}). Let us denote
this sign by $\varepsilon$. Also, let $\tau=\tau^\varphi(X,\omega)$, $\tau'=\tau^{\varphi'}(X',\omega')$ and $\tau''=\tau^{\varphi''}(X'',\omega'')$. By Lemma \ref{torsion} and the definition of the Conway function, the following
equalities hold in $Q(H)/H$:
\begin{eqnarray*}
-\varepsilon\,\nabla_L(t_1,\dots,t_n)&=&\varepsilon\,\tau(t_1^2,\dots,t^2_n)\\
&=&i'(\tau'({t'}^2,t_1^2,\dots,t^2_m))\,i''(\tau''({t''}^2,t_{m+1}^2,\dots,t^2_n))\\
&=&i'(-\nabla_{L'}(t',t_1,\dots,t_m))\,i''(-\nabla_{L''}(t'',t_{m+1},\dots,t_n))\\
&=&\nabla_{L'}(T',t_1,\dots,t_m)\,\nabla_{L''}(T'',t_{m+1},\dots,t_n),
\end{eqnarray*}
where $T'=i'(t')=t_{m+1}^{\ell''_{m+1}}\cdots t^{\ell''_n}_n$ and $T''=i''(t'')=t_1^{\ell'_1}\cdots t_m^{\ell'_m}$.
Therefore,
$$
\nabla_L(t_1,\dots,t_n)=-\varepsilon\,t_1^{\mu_1}\cdots
t_n^{\mu_n}\,\nabla_{L'}(T',t_1,\dots,t_m)\,\nabla_{L''}(T'',t_{m+1},\dots,t_n)
$$
in $Q(H)$, for some integers $\mu_1,\dots,\mu_n$. Now, the Conway function satisfies the symmetry fomula
$$
\nabla_L(t^{-1}_1,\dots,t^{-1}_n)=(-1)^n\,\nabla_L(t_1,\dots,t_n).
$$
Using this equation for $\nabla_L$, $\nabla_{L'}$ and $\nabla_{L''}$, it easily
follows that $\mu_1=\dots=\mu_n=0$. Therefore,
$$
\nabla_L(t_1,\dots,t_n)=-\varepsilon\,\nabla_{L'}(T',t_1,\dots,t_m)\nabla_{L''}(T'',t_{m+1},\dots,t_n)\in Q(H),
$$
where $\varepsilon$ is the sign of $\tau^{\varphi_T}(T,\omega_T)$. It remains to check that $\varepsilon=-1$. This can be
done by direct computation or by the following argument. Let $L'$ be the positive Hopf link in $S^3$, and let $L''$ be any
link such that $\nabla_{L''}\neq 0$. Clearly, the splice $L$ of $L'$ and $L''$ is equal to $L''$. Since $\nabla_{L'}=1$,
the equation above gives
$$
\nabla_{L''}(t_1,\dots,t_n)=-\varepsilon\,\nabla_{L''}(t_1,\dots,t_n).
$$
Since $\nabla_{L''}\neq 0$ and $\varepsilon$ does not depend on $L''$, we have $\varepsilon=-1$. This concludes the proof
of Theorem \ref{thm}.

\begin{acknowledgements}
The author thankfully acknowledges the Institut de Recherche Math\'ematique Avanc\'ee (Strasbourg) and the
Institut de Math\'ematiques de Bourgogne (Dijon) for hospitality. He also wishes to thank Vladimir Turaev and Mathieu
Baillif.
\end{acknowledgements}

\end{document}